\documentclass[graybox]{svmult}

\usepackage{amsfonts, amssymb, amsmath}
\usepackage{mathptmx}       
\usepackage{helvet}         
\usepackage{courier}        
\usepackage{type1cm}        
\usepackage{makeidx}         
\usepackage{graphicx}        
\usepackage{multicol}        
\usepackage[bottom]{footmisc}

\makeindex             

\begin{document}

\title*{Fitness,  Apprenticeship, and Polynomials}
\titlerunning{Fitness,  Apprenticeship, and Polynomials}
\author{Bernd Sturmfels}
\institute{
Bernd Sturmfels \at
Department of Mathematics, University of California, Berkeley, CA, 94720,
United
States of America  \ \ and \ \ 
Max-Planck Institute for Mathematics in the Sciences, 
Inselstra\ss e 22, 04103 Leipzig, Germany, \ \ \ 
\email{\texttt{bernd@berkeley.edu} \  or \
\texttt{bernd@mis.mpg.edu}}}
%
%
\maketitle


\abstract{This article discusses the design of the Apprenticeship 
Program at the Fields Institute, held 21 August--3 September 2016.
Six themes from combinatorial algebraic geometry
were selected for the two weeks: curves, surfaces, Grassmannians, convexity,
abelian combinatorics, parameters and moduli.
The activities were structured into
fitness, research and scholarship.
Combinatorics and concrete computations with polynomials 
(and theta functions) empowers young scholars in algebraic geometry,
and it helps them to connect with the historic roots of their field.
We illustrate our perspective for the threefold obtained by blowing up
six points in $\mathbb{P}^3$.
}

\section{Design}
\label{sec:1}

A thematic program on {\em Combinatorial Algebraic Geometry}
took place at the Fields Institute, Toronto, Canada, during the 
Fall Semester 2016. The program organizers were
David Cox, Megumi Harada, Diane Maclagan, Gregory Smith, and Ravi Vakil.

As part of this semester, the Clay Mathematics Institute funded the ``Apprenticeship Weeks'',
held 21 August--3 September 2016. This article discusses the
design and mathematical scope of this fortnight.
The structured activities took place in the mornings and afternoons
on Monday, Wednesday, and Friday, as well as the
mornings on Tuesday and Thursday.  The posted schedule was
identical for both weeks:

\smallskip

MWF 9:00--9:30: 	Introduction to today's theme

MWF 9:30--11:15:	Working on fitness problems

MWF 11:15--12:15:	Solutions to fitness problems

MWF 14:00--14:30:	Dividing into research teams

MWF 14:30--17:00:	Team work on projects

MWF 17:00--18:00:	Teams present findings

\smallskip

TuTh 9:00--12:00: Discussion of the scholarship theme

\medskip

\noindent
The term ``fitness'' is an allusion to physical exercise.
In order to improve physical fitness, many of us go to the gym.
A personal trainer can greatly enhance that experience.
The trainer develops your exercise plan and he pushes you
beyond previously perceived limits. The trainer makes you 
sweat a lot, he ensures that you use exercise equipment correctly,  and he
helps you to feel good about yourself afterwards.
In the context of team sports, the coach plays that role.
She works towards the fitness of the entire team,
where every player will contribute to the best of their abilities.

The six fitness sessions were designed to be as intense as those in sports.
Ten problems were posted for each session, and these were available
online two or three days in advance.
By design, these demanding problems were
open-ended and probed a different aspect of the theme. 
Section \ref{sec:3} of this article contains the complete list of problems,
along with a brief discussion and references that contain some solutions.

The ``apprentices'' were about $40$ early-career mathematicians, graduate students
and postdocs, coming from a wide range of backgrounds.
An essential feature of the Apprenticeship Weeks was the effort 
to build teams, and to promote collaboration as much as possible.
This created an amazing sense of community within the group.

At 9:00am on each Monday, Wednesday or Friday, a brief introduction was given to each 
fitness question. We formed ten teams to work on the problems.    At 11:15am
we got together again, and  one person from each team gave a brief presentation on
what had been discussed and discovered.
 Working on a challenging problem, with a group of new collaborators,
for less than two hours created a very intense and stimulating experience.
A  balanced selection process ensured that each participant 
had the opportunity to present for their team at least once.

At 2:00pm the entire group re-assembled and they discussed research-oriented
problems for the afternoons. This was conducted in the style of the
American Institute for Mathematics (AIM), whereby one of the participants
serves as the discussion leader, and only that person is allowed to touch the blackboard.
This led to an ample supply of excellent questions, some a direct continuation of
the morning fitness problems, and others only vaguely inspired by these.
Again, groups were formed for the afternoon, and they engaged in
learning and research. Computations and literature search played a big role,
and a lot of teaching went on in the groups.

Tuesday and Thursdays were discussion days. Here the aim was
to create a sense of scholarship among the participants.
The morning of these days involved
studying various software packages, classical research papers from the 19-th
and early 20-th centuries, and the diverse applications of combinatorial algebraic
geometry.  The prompts are given in Section \ref{sec:2}.
 The afternoons on discussion days were unstructured to allow the
participants time to ponder, probe, and write up their many new ideas.

\section{Scholarship Prompts}
\label{sec:2}

Combinatorial algebraic geometry is a field that, by design, straddles mathematical boundaries.
One aim is to study algebraic varieties with special combinatorial features. At its roots, this field
is about systems of polynomial equations in several variables, and about symmetries
and other special structures in their solution sets. 

Section 5 offers a concrete illustration of this perspective
 for a system of polynomials in $32$ variables.
 The objects of combinatorial algebraic geometry are
amenable to a wide range of software tools, which are now used widely among
the researchers. 

Another point we discussed is the connection to problems outside of pure mathematics.
A  new field, {\em Applied Algebraic Geometry}, has arisen in the past decade.
The techniques used there  often connect back 
to 19th and early 20th century work in algebraic geometry,
which is much more concrete and combinatorial than many recent
developments. And, even for her study of current abstract theories,
an apprentice may benefit from knowing the historic origins that have inspired the
development of algebraic geometry. Understanding these aspects, by getting
hands-on experiences and by studying original sources, was a focus in this
part of the program. 

\smallskip

In what follows we replicate the hand-outs for the four TuTh mornings.
The common thread can be summarized as: back to the roots.
These were given to the participants as prompts for explorations
and discussions. For several of the participants, it was their
first experience with software for algebraic geometry.
For others, it offered a first opportunity to 
 read an article that was published over 100 years ago.

\subsubsection*{Tuesday, August 23: Software}

Which software tools are  most useful for performing computations in \\
{\em Combinatorial Algebraic Geometry} ? Why?

\smallskip

Many of us are familiar with {\tt Macaulay2}.
Some of us are familiar with \\ {\tt Singular}.
What are your favorite packages within these systems?

\smallskip

Lots of math is supported by general-purpose computer algebra systems
such  as {\tt Sage}, {\tt Maple}, {\tt Mathematica}, or {\tt Magma}.
Do you use any of these regularly? For research or for teaching? 
How often and in which context?

\smallskip

Other packages that are useful for our community include
{\tt Bertini},
{\tt PHCpack}, 
{\tt 4ti2},
{\tt Polymake},
{\tt Normaliz},
{\tt GFan}.
What are these and what do they do?
Who developed them and why?

\smallskip

Does visualization matter in algebraic geometry? \\
Have you tried software like {\tt Surfex}?

\smallskip

Which software tool do you want to learn today?

\subsubsection*{Thursday, August 25: The 19th Century}

Algebraic Geometry has a deep and distinguished  history
that goes back hundreds of years. Combinatorics entered
the scene a bit more recently.

\medskip

Young scholars interested in algebraic geometry are strongly encouraged to familiarize
themselves with the literature from the 19th century.
Dig out papers from that period and {\em read them}!
Go for the original sources.  Some are in English.
Do not be afraid of languages like French, German, Italian.

Today we form groups. Each group will explore the life and work of
one mathematician, with focus on what he has done in
algebraic geometry.  Identify one key paper written by that author.
Then present your findings.

Here are some suggestions, listed alphabetically:
\begin{itemize}
\item Alexander von Brill 
\item Arthur Cayley
\item Michel Chasles 
\item Luigi Cremona 
\item Georges Halphen 
\item Otto Hesse 
\item Ernst Kummer 
\item Max Noether 
\item Julius Pl\"ucker 
\item Bernhard Riemann  
\item Friedrich Schottky 
\item Hermann Schubert 
\item Hieronymus Zeuthen 
\end{itemize}

\subsubsection*{Tuesday, August 30: Applications}

The recent years have seen a lot of interest in applications of
algebraic geometry, outside of core pure mathematics. An
influential
event was a research year 2006-07 at the IMA in Minneapolis. Following a suggestion by 
Doug Arnold (then IMA director and SIAM president), it led to the creation of the SIAM activity group
in Algebraic Geometry, and (ultimately) to the SIAM Journal on Applied Algebra and Geometry.
The reader is referred to these resources for further information.
These interactions with the
sciences and engineering have been greatly enhanced by the interplay with
Combinatorics and Computation  seen here at
the Fields Institute. However, the term ``Algebraic Geometry'' has to
be understood now in a broad sense.

\medskip

Today we form groups. Each group will get familiar with one
field of application, and they will select one paper in
Applied Algebraic Geometry that represent an interaction with that field.
Read your paper and then present your findings.
Here are some suggested fields, listed alphabetically:
\begin{itemize}
\item Approximation Theory 
\item Bayesian Statistics 
\item Chemical Reaction Networks 
\item Coding Theory  
\item Combinatorial Optimization 
\item Computer Vision  
\item Cryptography  
\item  Game Theory 
\item Geometric Modeling 
\item Machine Learning 
\item Maximum Likelihood Inference 
\item Neuroscience 
\item Phylogenetics 
\item Quantum Computing 
\item Semidefinite Programming  
\item Systems Biology 
\end{itemize}

\subsubsection*{Thursday, September 1: The Early 20th Century}

One week ago
 we examined the work of some algebraic geometers from
the 19th century. Today, we move on to the early 20th century,
to mathematics that was published prior to World War II.
You are encouraged to familiarize
yourselves with the literature from the period 1900-1939.
Dig out papers from that period and {\em read them}!
Go for the original sources.  Some are written in English.
Do not be afraid of languages like French, German, Italian, Russian.

Each group will explore the life and work of
one mathematician, with focus on what (s)he has done in 
algebraic geometry during that period. Identify one key 
paper written by that author. Then present your findings.

Here are some suggestions, listed alphabetically:
\begin{itemize}
\item  Eugenio Bertini 
\item Guido Castelnuovo 
\item Wei-Liang Chow 
\item Arthur B.~Coble 
\item Wolfgang Gr\"obner  
\item William V.D.~Hodge 
\item Wolfgang Krull 
\item  Solomon Lefschetz  
\item Frank Morley 
\item Francis S.~Macaulay 
\item Amalie Emmy Noether  
\item Ivan Georgievich Petrovsky  
\item Virginia Ragsdale 
\item Gaetano Scorza  
\item Francesco Severi 
\end{itemize}

\section{Fitness Prompts}
\label{sec:3}
  
 This section presents the six worksheets for the morning sessions
 on Mondays, Wednesdays and Fridays. These prompts inspired
 most of the articles in this volume. Specific pointers to dates
 refer to events that took place at the Fields Institute.
  The next section contains notes
 for each problem, offering references and solutions.
 
\subsubsection*{Monday, August 22: Curves}

\begin{enumerate}
\item Which genus can a smooth curve of degree $6$ in $\mathbb{P}^3$ have?
  Give examples.
\item Let $f(x) = (x-1)(x-2)(x-3)(x-6)(x-7)(x-8)$ and consider the genus $2$
  curve $y^2 = f(x)$.  Where is it in the moduli space $\mathcal{M}_2\,$?  Compute the Igusa
  invariants. Draw the Berkovich skeleton for the field of $5$-adic numbers.
\item The \emph{tact invariant} of two plane conics is the polynomial of
  bidegree $(6,6)$ in the $6+6$ coefficients which vanishes when the conics
  are tangent. Compute this invariant explicitly.  How many terms does it
  have?
\item \emph{Bring's curve} lives in a hyperplane in $\mathbb{P}^4$. It is defined by
  $ x_0^i + x_1^i + x_2^i + x_3^i + x_4^i = 0$ for $i = 1, 2, 3$.  What is its
  genus? Determine all tritangent planes of this curve.
\item Let $X$ be a curve of degree $d$ and genus $g$ in $\mathbb{P}^3$.  The
  Chow form of $X$ defines a hypersurface in the Grassmannian
  $\operatorname{Gr}(1,\mathbb{P}^3)$.  Points are lines that meet $X$.  Find
  the dimension and (bi)degree of its singular locus.
\item What are the equations of the secant varieties of \emph{elliptic normal
    curves}?
\item Let $X_P$ be the \emph{toric variety} defined by a $3$-dimensional
  lattice polytope, as in  Milena Hering's July 18-22 course.  Intersect $X_P$
  with two general hyperplanes to get a curve.  What is the degree and genus
  of that curve?
\item A 2009 article by Sean Keel and Jenia Tevelev presents \emph{Equations
    for $\overline{\mathcal{M}}_{0,n}$}.  Write these equations in \texttt{Macaulay2}
  format for $n=5$ and $n=6$.  Can you see the $\psi$-classes 
  (seen in Renzo Cavalieri's July 18-22 course)
     in these coordinates?
\item Review the statement of \emph{Torelli's Theorem} for genus $3$.  Using
  \texttt{Sage} or \texttt{Maple}, compute the $3 \times 3$ Riemann matrix of
  the Fermat quartic $\{x^4+y^4+z^4=0\}$.  How can you recover the curve from
  that matrix?
\item The moduli space $\mathcal{M}_7$ of genus $7$ curves has dimension $18$.
  What is the codimension of the locus of plane curves?  Hint: Singularities
  are allowed.
\end{enumerate}

\subsubsection*{Wednesday, August 24: Surfaces}

\begin{enumerate}
\item A nondegenerate surface in $\mathbb{P}^n$ has degree at least
  $n-1$. Prove this fact and determine all surfaces of degree $n-1$. Give
  their equations.
\item How many lines lie on a surface obtained by intersecting two quadratic
  hypersurfaces in $\mathbb{P}^4$?  Find an instance where all lines are
  defined over~$\mathbb{Q}$.
\item What is the maximum number of singular points on an irreducible quartic
  surface in $\mathbb{P}^3$?  Find a surface and compute its \emph{projective
    dual}.
\item Given a general surface of degree $d$ in $\mathbb{P}^3$, the set of its
  \emph{bitangent lines} is a surface in $\operatorname{Gr}(1,
  \mathbb{P}^3)$. Determine the cohomology class (or bidegree) of that surface.
\item Pick two random circles $C_1$ and $C_2$ in $\mathbb{R}^3$.  Compute
  their \emph{Minkowski sum} $C_1 + C_2$ and their \emph{Hadamard product}
  $C_1 \star C_2$.  Try other curves.
\item Let $X$ be the surface obtained by blowing up five general points in the plane.
  Compute the \emph{Cox ring} of $X$. Which of its ideals describe points on
  $X$?
\item The incidences among the $27$ lines on a cubic surface defines a
  $10$-regular graph.  Compute the complex of independent sets in this graph.
\item The Hilbert scheme of points on a smooth surface is smooth. Why?  How
  many torus-fixed points are there on the Hilbert scheme of $20$ points in
  $\mathbb{P}^2$? What can you say about the graph that connects them?
\item State the \emph{Hodge Index Theorem}. Verify this theorem for cubic
  surfaces in $\mathbb{P}^3$, by explicitly computing the matrix for the
  intersection pairing.
\item List the equations of one \emph{Enriques surface}. Verify its Hodge
  diamond.
\end{enumerate}

\subsubsection*{Friday, August 26: Grassmannians}

\begin{enumerate}
\item Find a point in $\operatorname{Gr}(3,6)$ with precisely $16$ non-zero
  Pl\"ucker coordinates.  As in June Huh's   July 18-22 course,
  determine the Chow  ring of its \emph{matroid}.
\item The coordinate ring of the Grassmannian $\operatorname{Gr}(3,6)$ is a
  \emph{cluster algebra} of finite type.  What are the cluster variables?
  List all the clusters.
\item Consider two general surfaces in $\mathbb{P}^3$ whose degrees are $d$
  and $e$ respectively.  How many lines in $\mathbb{P}^3$ are \emph{bitangent}
  to both surfaces?
\item The \emph{rotation group} $\operatorname{SO}(n)$ is an affine variety in
  the space of real $n \times n$-matrices. Can you find a formula for the
  degree of this variety?
\item The \emph{complete flag variety} for $\operatorname{GL}(4)$ is a
  six-dimensional subvariety of
  $\mathbb{P}^3 \times \mathbb{P}^5 \times \mathbb{P}^3$.  Compute its ideal
  and determine its tropicalization.
\item Classify all toric ideals that arises as initial ideals for the flag variety
  above. For each such toric degeneration, compute the \emph{Newton-Okounkov
    body}.
\item The Grassmannian $\operatorname{Gr}(4,7)$ has dimension $12$.  Four
  \emph{Schubert cycles} of codimension $3$ intersect in a finite number of
  points.  How large can that number be? Exhibit explicit cycles whose
  intersection is reduced.
\item The \emph{affine Grassmannian} and the \emph{Sato Grassmannian} are two
  infinite-dimensional versions of the Grassmannian. How are they related?
\item The coordinate ring of the Grassmannian $\operatorname{Gr}(2,7)$ is
  $\mathbb{Z}^7$-graded.  Determine the Hilbert series and the multidegree of
  $\operatorname{Gr}(2,7)$ for this grading.
\item The \emph{Lagrangian Grassmannian} parametrizes $n$-dimensional
  isotropic subspaces in $\mathbb{C}^{2n}$.  Find a Gr\"obner basis for its
  ideal. What is a `doset'?
\end{enumerate}

\subsubsection*{Monday, August 29: Convexity}

\begin{enumerate}
\item The set of nonnegative binary sextics is a closed full-dimensional
  convex cone in $\operatorname{Sym}_6(\mathbb{R}^2) \simeq
  \mathbb{R}^{7}$. Determine the face poset of this convex cone.
\item Consider \emph{smooth} projective toric fourfolds with eight invariant
  divisors.  What is the maximal number of torus-fixed points of any such
  variety?
\item Choose three general ellipsoids in $\mathbb{R}^3$ and compute the convex
  hull of their union. Which algebraic surfaces contribute to the boundary?
\item Explain how the Alexandrov-Fenchel Inequalities (for convex bodies) can
  be derived from the Hodge Index Theorem (for algebraic surfaces).
\item The blow-up of $\mathbb{P}^3$ at six general points is a threefold that
  contains $32$ special surfaces (exceptional classes).  What are these surfaces? Which triples
  intersect?  Hint: Find a $6$-dimensional polytope that describes the
  combinatorics.
\item Prove that every face of a spectrahedron is an exposed face.
\item How many combinatorial types of reflexive polytopes are there in
  dimension $3$?  In dimension $4$? Draw pictures of some extreme specimen.
\item A $4 \times 4$-matrix has six off-diagonal $2 \times 2$-minors.  Their
  binomial ideal in $12$ variables has a unique toric component.  Determine
  the f-vector of the polytope (with $12$ vertices) associated with this toric
  variety.
\item Consider the Pl\"ucker embedding of the real Grassmannian
  $\operatorname{Gr}(2,5)$ in the unit sphere in $\mathbb{R}^{10}$.  Describe
  its convex hull. Hint: Calibrations, Orbitopes.
\item Examine Minkowski sums of three tetrahedra in $\mathbb{R}^3$. What is
  the maximum number of vertices such a polytope can have? How to generalize?
\end{enumerate}

\subsubsection*{Wednesday, August 31:  Abelian Combinatorics}

\begin{enumerate}
\item The intersection of two quadratic surfaces in $\mathbb{P}^3$ is an
  \emph{elliptic curve}.  Explain its group structure in terms of geometric
  operations in $\mathbb{P}^3$.
\item A 2006 paper by Keiichi Gunji gives explicit equations for all
  \emph{abelian surfaces} in $\mathbb{P}^8$. Verify his equations in
  \texttt{Macaulay2}.  How to find the group law?
\item Experiment with Swierczewski's \texttt{Sage} code for the numerical
  evaluation of the \emph{Riemann theta function} $\theta(\tau;z)$.  Verify
  the functional equation.
\item \emph{Theta functions with characteristics}
  $\theta[\epsilon,\epsilon'](\tau;z)$ are indexed by two binary vectors
  $\epsilon, \epsilon' \in \{0,1\}^g$.  They are odd or even. How many each?
\item Fix the symplectic form
  $\langle x, y \rangle = x_1 y_4 + x_2 y_5 + x_3 y_6 + x_4 y_1 + x_5 y_2 + x_6
  y_3$
  on the $64$-element vector space $(\mathbb{F}_2)^6$.  Determine all
  isotropic subspaces.
\item Explain the combinatorics of the root system of type $\text{E}_7$.  How
  would you choose coordinates?  How many pairs of roots are orthogonal?
\item In 1879 Cayley published a paper in Crelle's journal titled $\,$
  \emph{Algorithms for ...}  What did he do? How does it relate the previous
  two exercises?
\item The \emph{regular matroid} $R_{10}$ defines a degeneration of abelian
  $5$-folds.  Describe its periodic tiling on $\mathbb{R}^5$ and secondary
  cone in the $2$-nd Voronoi decomposition.  Explain the application to Prym
  varieties due to Gwena.
\item Consider the Jacobian of the plane quartic curve defined over
  $\mathbb{Q}_2$ by
  \begin{align*}
    41x^4 &+ 1530x^3y + 3508x^3z + 1424x^2y^2 + 2490x^2yz \\
    & - 2274x^2z^2 + 470xy^3 + 680xy^2z  - 930xyz^2 + 772xz^3 \\
    & + 535y^4 - 350y^3z - 1960y^2z^2 - 3090yz^3 - 2047z^4 
  \end{align*} 
  Compute its limit in \emph{Alexeev's moduli space} for the $2$-adic
  valuation.
\item Let $\Theta$ be the \emph{theta divisor} on an abelian threefold $X$.
  Find $n = \dim H^0(X,k\Theta)$.  What is the smallest integer $k$ such that
  $k \Theta$ is very ample? Can you compute (in \texttt{Macaulay2}) the ideal of
  the corresponding embedding $X \hookrightarrow \mathbb{P}^{n-1}$?
\end{enumerate}

\subsubsection*{Friday, September 2: Parameters and Moduli}

\begin{enumerate}
\item Write down (in \texttt{Macaulay2} format) the two generators of the
  \emph{ring of invariants} for ternary cubics. For which plane cubics do both
  invariants vanish?
\item Fix a $\mathbb{Z}$-grading on the polynomial ring
  $S = \mathbb{C}[a,b,c,d]$ defined by $\deg(a) = 1$, $\deg(b) = 4$,
  $ \deg(c) = 5$, and $ \deg(d) = 9$.  Classify all homogeneous ideals $I$
  such that $S/I$ has Hilbert function identically equal to $1$.
\item Consider the Hilbert scheme of eight points in affine $4$-space
  $\mathbb{A}^4$.  Identify a point that is not in the main component. List
  its ideal generators.
\item Let $X$ be the set of all symmetric $4 \times 4$-matrices in
  $\mathbb{R}^{4 \times 4}$ that have an eigenvalue of multiplicity $\geq
  2$. Compute the $\mathbb{C}$-Zariski closure of $X$.
\item Which cubic surfaces in $\mathbb{P}^3$ are stable? Which ones are
  semi-stable?
\item In his second lecture on August 15, Valery Alexeev used six lines in
  $\mathbb{P}^2$ to construct a certain moduli space of K3 surfaces with $15$
  singular points. List the most degenerate points in the boundary of that
  space.
\item Find the most singular point on the Hilbert scheme of $16$ points in
  $\mathbb{A}^3$.
\item The polynomial ring $\mathbb{C}[x,y]$ is graded by the $2$-element group
  $\mathbb{Z}/2\mathbb{Z}$ where $\deg(x) = 1$ and $\deg(y) = 1$. Classify all
  Hilbert functions of homogeneous ideals.
\item Consider all threefolds obtained by blowing up six general points in
  $\mathbb{P}^3$.  Describe their Cox rings and Cox ideals.  How can you compactify
  this moduli space?
\item The moduli space of tropical curves of genus $5$ is a polyhedral space
  of dimension $12$. Determine the number of $i$-faces for
  $i = 0, 1, 2, \dotsc, 12$.
\end{enumerate}

\section{Notes, Solutions and References}
\label{sec:4}

Solutions to several of the sixty fitness problems can be found in the $16$
articles of this volume. The articles are listed as the first $16$ entries in our
References. They will be published in the order in which they are cited in this section.
In what follows we also offer references  for  other problems that did not lead
to articles in this book.

\subsubsection*{Notes on Curves}

\noindent 1. Castelnuovo classified  the
degree and genus pairs $(d,g)$ for all smooth curves in $\mathbb{P}^n$.
This was extended to characteristic $p$ by Ciliberto \cite{Cil}.
For $n=3, d =6$,  the  possible genera are $g=0,1,2,3,4$.
The {\tt Macaulay2} package {\tt RandomCurves} can compute examples.
The Hartshorne-Rao module \cite{Sch} plays a key role. \smallskip \\
\noindent 2. See Section 2 in the article by Bolognese, Brandt and Chua \cite{BBC}.
The approach using Igusa invariants was developed by
Helminck in \cite{Hel}. \smallskip \\
\noindent 3. The tact invariant has $3210$ terms, by \cite[Example 2.7]{StHurwitz}. \smallskip \\
\noindent 4. See Section 2.1 in the article by Harris and Len \cite{HL}.
The analogous problem for bitangents of plane quartics is discussed by
Chan and Jiradilok \cite{CJ}.
 \smallskip \\
\noindent 5. This is solved in the article by Kohn,  N{\o}dland and Tripoli \cite{KNT} \smallskip \\
\noindent 6. Following Fisher \cite{Fis},
elliptic normal curves are defined by the
$4 \times 4$-subpfaffians of the  Klein matrix,
and their secant varieties are defined by its larger subpfaffians.
 \smallskip \\
\noindent 7. The degree of a projective toric variety $X_P$ is the volume of its lattice polytope $P$.
The genus of a complete intersection in $X_P$
was derived by Khovanskii in 1978. We recommend the
tropical perspective offered by Steffens and Theobald in \cite[\S 4.1]{ST}.
 \smallskip \\
\noindent 8. See the article by Monin and Rana \cite{MR} for a solution up to $n=6$. \smallskip \\
\noindent 9. See \cite{DH} for how to compute the forward direction of the Torelli map 
of an arbitrary plane curve. For computing the backward direction in genus $3$ see \cite[\S 5.2]{SD}.
 \smallskip \\
\noindent 10. 
  Trinodal sextics form a $16$-dimensional family; their codimension in $\mathcal{M}_7$ is two.
  This is a result  due to Severi, derived by
  Castryck and Voight in \cite[Theorem 2.1]{CV}.

\subsubsection*{Notes on Surfaces}

\noindent 1. This was solved by Del Pezzo in 1886. 
Eisenbud and Harris \cite{EH} give a beautiful introduction to the
theory of {\em varieties of minimal degree}, including their equations.
 \smallskip \\
\noindent 2. This is a {\em del Pezzo surface} of degree $4$. It has $16$ lines.
To make them rational, map $\mathbb{P}^2$ into $\mathbb{P}^3$ via a $\mathbb{Q}$-basis for the
cubics that vanish at five rational points in $\mathbb{P}^2$.
 \smallskip \\
\noindent 3. The winner, with $16$ singular points, is the {\em Kummer surface} \cite{Hudson}. 
 It is self-dual. \smallskip  \\
\noindent 4. This is solved in the article by Kohn,  N{\o}dland and Tripoli \cite{KNT}.  \smallskip \\
\noindent 5. See Section 5 in the article by  Friedenberg, Oneto and Williams \cite{FOW}. \smallskip \\
\noindent 6. This is the del Pezzo surface in Problem 2. Its Cox ring is a polynomial ring
in $16$ variables modulo an ideal generated by $20$ quadrics.
Ideal generators that are universal over the base $\overline{\mathcal{M}_{0,5}}$ 
are listed in \cite[Proposition 2.1]{RSS}.
Ideals of  points on the surface are  torus translates of the
toric ideal of the $5$-dimensional demicube $D_5$.
For six points in $\mathbb{P}^2$ we refer to
Bernal, Corey,  Donten-Bury, Fujita and Merz \cite{BCDFM}.   \smallskip \\
\noindent 7. This is the clique complex of the {\em Schl\"afli graph}. 
The f-vector of this simplicial complex is  $(27, 216, 720, 1080, 648, 72)$.
The Schl\"afli graph is the edge graph of the {\em $E_6$-polytope}, denoted $2_{21}$,
which is a cross section of the {\em Mori cone} of  the surface.
 \smallskip \\
\noindent 8. The torus-fixed points on ${\rm Hilb}_{20}(\mathbb{P}^2)$ are indexed by ordered triples
of partitions $(\lambda_1,\lambda_2,\lambda_3)$ with $|\lambda_1|+|\lambda_2| +|\lambda_3|=20$.
The number of such triples equals $341,649$. The graph that connects them is a variant
of the graph for the Hilbert scheme of points in the affine plane.
The latter was studied by Hering and Maclagan in~\cite{HM}.
  \smallskip \\
\noindent 9. The signature of the intersection pairing is $(1,r-1)$
where $r$ is the rank of the Picard group. This is $r=7$ for the cubic surface.
From the analysis in Problem 7, we can get various symmetric matrices that 
represent the intersection pairing.
 \smallskip \\
\noindent 10. See the article by Bolognese, Harris and Jelisiejew \cite{BHJ}.  \smallskip 

\subsubsection*{Notes on Grassmannians}

\noindent 1. See the article by  Wiltshire-Gordon, Woo and Zajackowska \cite{WWZ}.  \smallskip \\
\noindent 2. In addition to the $20$ Pl\"ucker coordinates $p_{ijk}$,
one needs two more functions, namely
$p_{123} p_{456} - p_{124} p_{356}$ and
$p_{234} p_{561} - p_{235} p_{461}$. The six boundary Pl\"ucker coordinates
$p_{123}, p_{234}, p_{345}, p_{456}, p_{561}, p_{612}$ are frozen. The
other $16$ coordinates are the {\em cluster variables} for ${\rm Gr}(3,6)$.
This was derived by Scott in \cite[Theorem 6]{Sco}.
 \smallskip \\
\noindent 3. This is worked out in the article by Kohn,  N{\o}dland and Tripoli \cite{KNT}. \smallskip \\
\noindent 4. This is the main result of Brandt, Bruce, Brysiewicz, Krone and Robeva \cite{BBBKR}.
 \smallskip \\
\noindent 5. See 
the article by  Bossinger, Lamboglia, Mincheva and Mohammadi \cite{BLMM}.  \smallskip \\
\noindent 6. See 
the article by  Bossinger, Lamboglia, Mincheva, Mohammadi \cite{BLMM}.  \smallskip \\
\noindent 7. The maximum number is $8$. This is obtained by taking the
partition $(2,1)$ four times. For this  problem, and many other Schubert problems,
 instances exist where all solutions are real. 
See the works of Sottile, specifically \cite[Theorem 3.9 (iv)]{Sot}.
\smallskip \\
\noindent 8. The Sato Grassmannian is more general than the affine Grassmannian.
These are studied, respectively, in
{\em integrable systems} and in
{\em geometric representation theory}.
 \smallskip \\
 \noindent 9. A formula for the $\mathbb{Z}^n$-graded Hilbert series of
 ${\rm Gr}(2,n)$ is given by Witaszek \cite[\S 3.3]{Wit}. 
 For an introduction to multidegrees see  \cite[\S 8.5]{CCA}.
  Try the {\tt Macaulay2} commands {\tt Grassmannian} and {\tt multidegree}.
Escobar and Knutson \cite{EK} determine the multidegree of a
variety that is important in computer vision.
  \smallskip \\
\noindent 10. The coordinate ring of the Lagrangian Grassmannian is an algebra
with straightening law over a {\em doset}. This stands for double poset.
See the exposition in \cite[\S 3]{Ruf}.

\smallskip

\subsubsection*{Notes on Convexity}

\noindent 1. The face lattice of the cone of 
non-negative binary forms of degree $d$ is described
in Barvinok's textbook \cite[\S II.11]{Bar}. In more variables this is much more difficult.
   \smallskip \\
\noindent 2. This seems to be an open problem. For seven invariant divisors,
this was resolved by Gretenkort {\it et al.} \cite{GKS}.
Note the conjecture stated in the last line  of that paper.
  \smallskip \\
\noindent 3. We refer to Nash, Pir, Sottile and Ying \cite{NPSY} and to the
{\tt youtube} video {\em The Convex Hull of Ellipsoids} by 
Nicola Geismann, 
Michael Hemmer,
and Elmar Sch\"omer.
 \smallskip \\
\noindent 4. We refer to Ewald's textbook, specifically \cite[\S IV.5 and \S VII.6]{Ewa}. \smallskip \\
\noindent 5. 
The relevant polytope is the $6$-dimensional demicube;
its $32$ vertices correspond to the $32$ special divisors.
See the notes for Problem 9 in Parameters and Moduli.
 \smallskip \\
\noindent 6. 
This was first proved by Ramana and Goldman in \cite[Corollary 1]{RG}.
 \smallskip \\
\noindent 7. 
Kreuzer and Skarke \cite{KS} classified such reflexive polytopes
up to lattice isomorphism. There are $4319$ in dimension $3$, and there are
$473800776$ in dimension~$4$.  Lars Kastner classified the 
list of $4319$ into combinatorial types.
He found that  there are $558$ combinatorial types of reflexive $3$-polytopes.
They have up to $14$ vertices.
 \smallskip \\
\noindent 8. This $6$-dimensional polytope is obtained from the direct product of two 
identical regular
tetrahedra by removing the four pairs of corresponding vertices. It is the
convex hull of the points $\,e_i \oplus e_j \,$ in $\,\mathbb{R}^4 \oplus \mathbb{R}^4$
  where $i,j \in \{1,2,3,4\}$ with $i \not= j$. Using the software 
  {\tt Polymake}, we find its f-vector to be  $(12,54, 110, 108, 52, 12)$.
  \smallskip \\
\noindent 9. 
The faces of the Grassmann orbitopes ${\rm conv}({\rm Gr}(2,n))$ 
for $n \geq 5$ are described in \cite[Theorem 7.3]{SSS}. 
It is best to start with the easier case $n=4$ 
in \cite[Example 7.1]{SSS}.
  \smallskip \\
\noindent 10. The maximum number of vertices is $38$, by the formula
of Karavelas {\it et al.}~in
\cite[\S 6.1, equation (49)]{KKT}.
A definitive solution to the problem of
characterizing face numbers of Minkowski sums of polytopes
was given by Adiprasito and Sanyal~\cite{AS}.
\smallskip

\subsubsection*{Notes on Abelian Combinatorics}
\noindent 1. A beautiful solution was written up by Qiaochu Yuan when he was a high school student;
see \cite{Yuan}. The idea is to simultaneously diagonalize
the two quadrics, then project their intersection curve into the plane, thereby
obtaining an {\em Edwards curve}.  \smallskip \\
\noindent 2.  This  is a system of  $9$ quadrics and $3$ cubics, derived
 from Coble's cubic as in \cite[Theorem 3.2]{RSamS}.
Using theta functions as in
\cite[Lemma 3.3]{RSamS}, one gets the group law.
\smallskip \\
\noindent 3. 
See \cite{SD} and compare with Problem 9 in Curves.
  \smallskip \\
\noindent 4. For the $2^{2g}$ pairs $(\epsilon,\epsilon')$,
we check whether $\epsilon \cdot \epsilon'$ is even or odd.
There are $2^{g-1}(2^g+1)$ even theta characteristics
and $2^{g-1}(2^g-1)$ odd theta characteristics.   \smallskip \\
\noindent 5. The number of isotropic subspaces of $(\mathbb{F}_2)^6$ is $63$ of dimension $1$, it is
$315$ in dimension $2$, and it is $135$ in dimension $3$.  
The latter are the Lagrangians \cite[\S 6]{RSSS}. \smallskip \\
\noindent 6.
The root system of type $E_7$ has $63$ positive roots. 
They are discussed in \cite[\S 6]{RSSS}.
\smallskip \\
\noindent 7. Cayley gives a bijection between the $63$ positive roots of $E_7$ with
the $63$ non-zero vectors in
$(\mathbb{F}_2)^6$. Two roots have inner product zero if
and only if the corresponding vectors in  $(\mathbb{F}_2)^6 \backslash \{ 0 \}$
are orthogonal in the setting of Problem 5. See \cite[Table 1]{RSSS}.
  \smallskip \\
\noindent 8. This refers to Gwena's article \cite{Gwe}. Since the matroid $R_{10}$
is not co-graphic, the corresponding tropical abelian varieties 
are not in the Schottky locus of Jacobians.
  \smallskip \\
\noindent 9. 
This fitness problem is solved in the article by Bolognese, Brandt and Chua \cite{BBC}
Chan and Jiradilok \cite{CJ} study an important special family of plane quartics.
   \smallskip \\
\noindent 10. The divisor $k \Theta$ is very ample
for $k = 3$.
This embeds any abelian threefold into $\mathbb{P}^{26}$.
For products of three cubic curves, each in $\mathbb{P}^2$, this gives
the Segre embedding.

\subsubsection*{Notes on Parameters and Moduli}

\noindent 1. The solution can be found, for instance,  on the website
$$ \hbox{\tt http://math.stanford.edu/$\sim$notzeb/aronhold.html}$$
The two generators have degree $4$ and $6$. The quartic invariant is known as
the {\em Aronhold invariant} and it vanishes when the ternary cubic is a sum of three cubes
of linear forms. Both invariants vanish when the cubic curve has a cusp.  \smallskip \\
\noindent 2. This refers to extra irreducible components in toric Hilbert schemes \cite{PS}.
These schemes were first introduced by  Arnold \cite{Arn},
who coined the term {\em A-graded algebras}.
Theorem 10.4 in \cite{GBCP} established the existence of an extra component for
$A = (1347)$.  We ask to verify
 the second entry in Table 10-1 on page 88 of \cite{GBCP}.
 \smallskip \\
\noindent 3. Cartwright {\it et al.}~\cite{CEVV} showed
that the Hilbert scheme of eight points in $\mathbb{A}^4$ has two
irreducible components. An explicit point in the non-smoothable component~is 
given in the article by Douvropoulos, Jelisiejew, N{\o}dland and Teitler \cite{DJNT}. 
 \smallskip \\
\noindent 4. At first, it is surprising that $X$  has codimension $2$.
The point is that we work over the real numbers
$\mathbb{R}$. The analogous set over $\mathbb{C}$
is the hypersurface of  a  sum-of-squares polynomial.
 The $\mathbb{C}$-Zariski closure of
$X$ is a nice variety of codimension~$2$. 
The defining ideal and its Hilbert-Burch resolution are explained in
\cite[\S 7.1]{CBMS}.  \smallskip \\
\noindent 5. This is an exercise in Geometric Invariant Theory \cite{MFK}.
 A cubic surface is stable if and only if it has at most ordinary double points
 (A1 singularities). For semi-stable surfaces,  A2 singularities are allowed.
 For an exposition see \cite[Theorem 3.6]{Rei};
 this is E.~Reinecke's Bachelor thesis,  written under the supervision of
 D.~Huybrechts. \smallskip \\
 \noindent 6. This is the moduli space of stable hyperplane arrangements \cite{Ale},
 here for the case of six lines in $\mathbb{P}^2$. The precise space depends on a choice of
  parameters \cite[\S 5.7]{Ale}. For some natural parameters, this is the tropical
 compactification associated with the tropical Grassmannian ${\rm Gr}(3,6)$,
 so the most degenerate points correspond to the seven generic types of
 tropical planes in $5$-space, shown in  \cite[Figure 5.4.1]{MS}.  \smallskip \\ 
\noindent 7. See \cite[Theorem 2.3]{four}. \smallskip \\
\noindent 8.  For each partition, representing a monomial ideal in $\mathbb{C}[x,y]$,
we count  the odd and even boxes in its Young diagram. The resulting Hilbert functions
$h:\mathbb{Z}/2\mathbb{Z} \rightarrow \mathbb{N}$ are
$(h({\rm even}), h({\rm odd})) = (k^2 + m, k(k+1) + m) $ or $ ((k+1)^2 + m, k(k+1) + m)$,
where $k,m \in \mathbb{N}$.
This was contributed by Dori Bejleri.
 For more details see \cite[\S 1.3]{BZ}.
\smallskip \\
\noindent 9. The blow-up of $\mathbb{P}^{n-3}$ at $n$ points is a
Mori dream space. Its Cox ring has $2^{n-1}$ generators, constructed
explicitly by Castravet and Tevelev in \cite{CT}. These form
a  Khovanskii basis \cite{KM}, by \cite[Theorem~7.10]{SX}. The Cox ideal is studied in \cite{SV}.
Each point on its variety represents a
rank two stable quasiparabolic vector bundle on $\mathbb{P}^1$ with 
$n$ marked points. The relevant moduli space is $\overline{\mathcal{M}_{0,n}}$.
 \smallskip \\
\noindent 10. 
The moduli space of tropical curves of genus $5$ serves as
 the first example in the article by Lin and Ulirsch \cite{LU}.
The article by Kastner, Shaw and Winz \cite{KSW} discusses
state-of-the-art software tools
for computing with such polyhedral spaces.

\section{Polynomials}
\label{sec:5}

The author of this article holds the firm belief that algebraic 
geometry concerns the study of solution sets to systems of polynomial equations.
Historically, geometers explored curves and surfaces that are zero sets of
polynomials. It is the insights gained from these basic figures that have led, over the course
of centuries,  to the profound depth and remarkable breadth of contemporary algebraic geometry.
However, many of the current theories are now far removed from explicit varieties,
and polynomials are nowhere in sight.  What we are advocating is
for algebraic geometry to take an outward-looking perspective. Our readers
 should be aware of the wealth of applications in the sciences and engineering, and be open to
a ``back to the basics'' approach in both teaching and scholarship. From this
perspective, the interaction with combinatorics can be particularly valuable.
Indeed, combinatorics is known to some as the ``nanotechnology of mathematics''.
It is all about explicit objects, those that can be counted, enumerated, and dissected
with laser precision. And, these objects include 
some beautiful polynomials and the ideals they generate.

The following example serves as an illustration. We work in a polynomial ring
$\mathbb{Q}[p]$ in $32$ variables, one for each subset of $\{1,2,3,4,5,6\}$ whose
cardinality is odd:
$$ p_1,p_2,\ldots,p_6, \,p_{123}, p_{124}, p_{125}, \ldots, p_{356}, p_{456},\,
p_{12345}, p_{12346}, \ldots, p_{23456} . $$
The polynomial ring $\mathbb{Q}[p]$ is $\mathbb{Z}^7$-graded by setting
${\rm degree}(p_\sigma) = e_0 + \sum_{i \in \sigma} e_i$, where
$e_0,e_1,\ldots,e_6$ is the standard basis of $\mathbb{Z}^7$.
Let $X$ be a $5 \times 6$-matrix of variables, and let $I$ be the kernel of
the ring map $ \mathbb{Q}[p]  \rightarrow \mathbb{Q}[X]$ that takes the variables
$p_\sigma$ to the determinant of the submatrix of $X$ with column indices $\sigma$
and row indices $1,2,\ldots,|\sigma|$.

The ideal $I$ is prime and $\mathbb{Z}^7$-graded.
 It has multiple geometric interpretations.
First of all, it describes the partial flag variety of points in $2$-planes in hyperplanes in
$\mathbb{P}^5$. This flag variety lives in
$\mathbb{P}^5 \times \mathbb{P}^{19} \times \mathbb{P}^5$,
thanks to the Pl\"ucker embedding.
Its projection into the  factor $\mathbb{P}^{19}$ is the 
Grassmannian ${\rm Gr}(3,6)$ of $2$-planes in $\mathbb{P}^5$.
Flag varieties are studied by
Bossinger, Lamboglia, Mincheva and Mohammadi in \cite{BLMM}.

But, let the allure of polynomials now speak for itself.
Our ideal $I$ has $66$ minimal quadratic generators. Sixty 
generators are unique up to scaling in their degree:
$$ 
\begin{matrix} {\rm degree}  & \qquad &  \hbox{\rm ideal generator} \\
(2,\,0,0,1,1,1,1)  & & p_3 p_{456} - p_4 p_{356} + p_5 p_{346} - p_6 p_{345}  \\
(2, \, 0,1,0,1,1,1) & & p_2 p_{456} - p_4 p_{256} + p_5 p_{246} - p_6 p_{245}  \\
\cdots  & & \cdots \quad \cdots \quad \cdots \\
(2, \, 1,1,1,1,0,0) & &  p_1 p_{234} - p_2 p_{134} + p_3 p_{124} - p_4 p_{123}  \smallskip \\
(2, \,0,1,1,1,1,2)  & &  p_{256} p_{346}-p_{246} p_{356}+p_{236} p_{456} \\
\cdots  & & \cdots \quad \cdots \quad \cdots \\
(2,\,2,1,1,1,1,0) & & p_{125} p_{134} -p_{124} p_{135} +p_{123} p_{145} \smallskip \\
(2,\,1,1,1,1,2,2) & &   p_{156} p_{23456}- p_{256}p_{13456}+ p_{356}p_{12456}- p_{456}p_{12356} \\
\cdots  & & \cdots \quad \cdots \quad \cdots \\
(2,\,2,2,1,1,1,1) & &   p_{123} p_{12456}- p_{124} p_{12356} +p_{125}p_{12346} - p_{126}p_{12345}
\end{matrix}
$$
The other six minimal generators live in degree $(2, \,1,1,1,1,1,1)$. 
These are the $4$-term Grassmann-P\"ucker relations, like
$p_{126} p_{345} -p_{125} p_{346}+p_{124} p_{356} -p_{123} p_{456}$.

Here is an alternate interpretation of the ideal $I$.
It defines a variety of dimension $15 = \binom{6}{2}$ in
$\mathbb{P}^{31}$ known as the {\em spinor variety}.
In this guise,  $I$ encodes the algebraic
relations among the principal subpfaffians
of a skew-symmetric $6 \times 6$-matrix. Such
subpfaffians are indexed with the subsets
of $\{1,2,3,4,5,6\}$ of even cardinality.
The trick is to fix a natural bijection
between even and odd subsets.
This variety is similar to the
{\em Lagrangian Grassmannian} seen
in fitness problem \# 10 on Grassmannians.

At this point, readers who like combinatorics
and computations may study $I$. Can you compute the
tropical variety of $I$? Which of its maximal cones are
{\em prime} in the sense of Kaveh and Manon \cite[Theorem 1]{KM}?
These determine {\em Khovanskii bases}
for $\mathbb{Q}[p]/I$ and hence toric degenerations of the spinor variety in $\mathbb{P}^{31}$.
Their  combinatorics is recorded 
in a list of {\em Newton-Okounkov polytopes} with $32$ vertices.

Each of these polytopes comes with a linear projection  to
the $6$-dimensional demicube, which is the convex hull in $\mathbb{R}^7$
of the $32$ points ${\rm deg}(p_\sigma)$. We saw this demicube
in fitness problem \# 5 on Convexity, whose theme we turn to shortly.

It is the author's opinion that
Khovanskii bases deserve more attention than the Newton-Okounkov bodies they give rise to.
The former are the algebraic manifestation of a toric degeneration.
These must be computed and verified. Looking at a Khovanskii basis
through the lens of convexity reveals the Newton-Okounkov body.

We now come to a third, and even more interesting, 
geometric interpretation of our $66$ polynomials.
It has to do with {\em Cox rings}, and their Khovanskii bases,
similar to those in the article by
Bernal, Corey,  Donten-Bury, Fujita and Merz.
We begin by replacing the generic $5 \times 6$-matrix $X$ by one
that has the special form in \cite[(1.2)]{CT}:
$$ X \quad = \quad \begin{pmatrix} 
    u_1^2 x_1 &   u_2^2 x_2 &   u_3^2 x_3 &   u_4^2 x_4 &   u_5^2 x_5 &   u_6^2 x_6 \\
       u_1 y_1 &       u_2 y_2 &     u_3 y_3 &     u_4 y_4 &     u_5 y_5 &     u_6 y_6  \\
u_1 v_1 x_1 & u_2 v_2 x_2 & u_3 v_3 x_3 & u_4 v_4 x_4 & u_5 v_5 x_5 & u_6 v_6 x_6 \\
       v_1 y_1 &       v_2 y_2 &       v_3 y_3 &       v_4 y_4 &       v_5 y_5 &       v_6 y_6  \\
   v_1^2 x_1 &    v_2^2 x_2 &    v_3^2 x_3 &    v_4^2 x_4 &    v_5^2 x_5 &    v_6^2 x_6 
   \end{pmatrix}.
   $$
   Now, the polynomial ring $\mathbb{Q}[X]$ gets replaced by 
   $\Bbbk[x_1,x_2,\ldots,x_6,y_1,y_2,\ldots,y_6]$
   where $\Bbbk$ is the field extension of $\mathbb{Q}$ generated by the entries of a 
   $2 \times 6$-matrix of scalars:
   \begin{equation}
   \label{eq:Umatrix}
   U \quad = \quad \begin{pmatrix}
   u_1 & u_2 & u_3 & u_4 & u_5 & u_6 \\  v_1 & v_2 & v_3 & v_4 & v_5 & v_6 
 \end{pmatrix}.
 \end{equation}
 We assume that the $2 \times 2$-minors of $U$ are non-zero. Let
 $J$ denote the kernel of the odd-minors map $\Bbbk[p] \rightarrow \Bbbk[X]$ as before.
The ideal $J$ is also $\mathbb{Z}^7$-graded and it strictly contains the ideal $I$.
 Castravet and Tevelev \cite[Theorem 1.1]{CT} proved that
$\Bbbk[p]/J$ is the {\em Cox ring} of the blow-up of $\mathbb{P}_\Bbbk^3$ at six points.
These  points are Gale dual to $U$.
We refer to $J$ as the {\em Cox ideal} of that rational threefold whose
Picard group $\mathbb{Z}^7$ furnishes the grading. The affine variety in 
$\mathbb{A}_\Bbbk^{32}$ defined by $J$ is 
$10$-dimensional (it is the {\em universal torsor}). Quotienting 
 by a $7$-dimensional torus action yields our threefold.
 The same story for blowing up five points in $\mathbb{P}_\Bbbk^2$
 is problem \# 6 on Surfaces.
 
 In  \cite{SV}  we  construct the Cox ideal by duplicating  the
 ideal of the spinor variety:
\begin{equation}
\label{eq:J=I+}
 J \,\,\,= \,\,\, I \,+\, {\bf u}* I .
\end{equation}
Here ${\bf u}$ is a vector in $(K^*)^{32}$ that is derived
from $U$. The ideal ${\bf u}*I$ is obtained from~$I$ by scaling the
variables $f_\sigma$ with the coordinates of ${\bf u}$. In particular,
the Cox ideal~$J$ is minimally generated by $132 $ quadrics.
Now, there are two generators in each of the sixty $\mathbb{Z}^7$-degrees
in our table, and there are $12$ generators in degree $(2,\,1,1,1,1,1,1)$.

Following \cite[Example~7.6]{SX}, we fix the rational function field $\Bbbk = \mathbb{Q}(t)$ and set
$$ U \,\, = \,\,  \begin{pmatrix}\, 1 &\, t &\, t^2\, &\, t^3 \,& \, t^4\, &\, t^5 \,\\
\,t^5\, & \,t^4 \,& \,t^3 \,& \,t^2\, & t &1  \end{pmatrix}. $$
The ring map $\Bbbk[p] \rightarrow \Bbbk[X]$ now maps the variables $p_\sigma$ like this:
$$ 
\begin{matrix}
p_1 \, \mapsto \,  \underline{x_1} \quad \qquad \qquad  \qquad \qquad 
 \qquad \qquad \qquad \qquad  \qquad \qquad \qquad  \qquad\, \,\,\,\\
p_{123} \, \mapsto \,
\underline{x_1 y_2 x_3} \,t^6
-(x_1 x_2 y_3 +y_1 x_2 x_3)\, t^7
+(y_1 x_2 x_3 + x_1 x_2 y_3)\, t^9
- x_1 y_2 x_3 \, t^{10} \quad \\
p_{12345} \,\mapsto \,
\underline{x_1 y_2 x_3 y_4 x_5}\, t^{10}  
-(y_1 x_2 x_3 y_4 x_5 + x_1 y_2 x_3 x_4 y_5 + \cdots + x_1 x_2 y_3 y_4 x_5)\, t^{11} + 
\cdots \\
\end{matrix}
$$
Here is a typical example of a $\mathbb{Z}^7$-degree with two minimal ideal generators:
$$
\begin{matrix}
(2, \, 1,1,1,1,0,0) & &  p_1 p_{234} - p_2 p_{134} + p_3 p_{124} - p_4 p_{123}  \smallskip \\
(2, \, 1,1,1,1,0,0) & & \,\,\,\,t^4 p_1 p_{234} - t^2 p_2  p_{134} + t^2 p_1 p_{234}  + p_4 p_{123}
\end{matrix}
$$
The algebra generators $p_\sigma$ form a Khovanskii basis for $\Bbbk[p]/J$
with respect to the $t$-adic valuation.
The toric algebra resulting from this flat family is generated
by the underlined monomials. Its toric ideal ${\rm in}(J)$ is generated
by $132$ binomial quadrics:
$$ 
\begin{matrix} {\rm degree}  & \qquad &  \hbox{\rm pair of binomial generators for ${\rm in}(J)$} \\
(2,\,0,0,1,1,1,1)  & & p_3 p_{456} - p_4 p_{356} \qquad \quad p_5 p_{346} - p_6 p_{345} \\
(2, \, 0,1,0,1,1,1) & & p_2 p_{456} - p_4 p_{256} \qquad \quad p_5 p_{246} - p_6 p_{245} \\
\cdots  & & \cdots \qquad \qquad \qquad \qquad \cdots \\
(2, \, 1,1,1,1,0,0) & &  p_1 p_{234} - p_2 p_{134} \qquad  \quad p_3 p_{124} - p_4 p_{123} \\
(2, \,0,1,1,1,1,2)  & &   p_6 p_{23456} - p_{236} p_{456} \qquad p_{246} p_{356}- p_{256} p_{346} \\  
\cdots  & & \cdots \qquad \qquad \qquad \qquad \cdots \\
(2,\,2,1,1,1,1,0) & & 
p_1 p_{12345} - p_{123} p_{145} \qquad p_{124} p_{135} - p_{125} p_{134}
\smallskip \\
(2,\,1,1,1,1,2,2) & & 
p_{156} p_{23456} - p_{256} p_{13456}   \quad
p_{356} p_{12456} - p_{456} p_{12356}  \smallskip \\
\cdots  & & \cdots \qquad \qquad \qquad \qquad \cdots \\
(2,\,2,2,1,1,1,1) & &
p_{123} p_{12456} - p_{124} p_{12356} \quad p_{125} p_{12346} - p_{126} p_{12345}
\end{matrix}
$$
These $132$ binomials define a toric variety that 
is a degeneration of our universal torsor.
The ideal ${\rm in}(J)$ is relevant in both
biology and physics. It represents the {\em Jukes-Cantor model} in phylogenetics \cite{SS}
and the {\em Wess-Zumino-Witten model} in conformal field theory \cite{Man}.
Beautiful polynomials can bring the sciences together.

Let us turn to another fitness problem.
The past three pages offered a capoeira approach
to \# 9 in Parameters and Moduli.
The compactification is that given by the tropical variety of the universal Cox ideal,
to be computed 
as in \cite{RSamS, RSS}.
The base space is $\overline{\mathcal{M}_{0,6}}$,
with points represented by $2 \times 6$-matrices $U$ as in (\ref{eq:Umatrix}).
We encountered several themes that are featured in other
articles in this book: flag varieties, Grassmannians, $\mathbb{Z}^n$-gradings, Cox rings,
Khovanskii bases, and toric ideals.
The connection to spinor varieties was developed in the article \cite{SV} with Mauricio Velasco.
The formula (\ref{eq:J=I+}) is derived in \cite[Theorem 7.4]{SV} for
the blow-up of $\mathbb{P}^{n-3}$ at $n$ points when $n \leq 8$. It is still a conjecture
for $n \geq 9$.
On your trail towards solving such open problems,  fill your backpack with polynomials.
They will guide you.

\begin{acknowledgement}
This article benefited greatly from comments by
Lara Bossinger, Fatemeh Mohammadi, 
Emre Sert\"oz, Mauricio Velasco
and an anonymous referee.
The apprenticeship program at the Fields Institute was
supported by the Clay Mathematics Institute. 
The author also acknowledges partial support from the Einstein Foundation Berlin,
MPI Leipzig, and the US National Science Foundation (DMS-1419018).
\end{acknowledgement}

\end{document}